# Cauchy-like Integral Theorems for Quaternion and Biquaternion Functions

R.A.W. Bradford

University of Bristol, Faculty of Engineering

**Abstract**

It is well known that there is an integral theorem for quaternion-valued functions analogous to Cauchy's Theorem for complex-valued functions, namely Fueter's Theorem. The class of quaternionic functions for which this applies are generally referred to as "regular functions", and these provide the most productive means of generalising the class of holomorphic complex functions. This paper derives a second integral theorem, also analogous to Cauchy's Theorem, and which is believed to be quite distinct from that of Fueter, despite appearances. The paper takes the opportunity to present the basis of the derivation of both theorems, and also their extension to the associated classes of right-regular and conjugate regular functions. Both theorems can also be extended into the biquaternionic domain in which the four quaternion coordinates may be complex valued. This is of interest in physics as the Hermitian biquaternions have a natural norm which is Minkowskian and provide an elegant formalism for Lorentz transformations.

MSC (2020) class: 30E20 (primary), 46S05 (secondary)

## Introduction

It is a curious fact that neither Hamilton nor his immediate followers in the mid-nineteenth century developed the analysis of quaternionic functions in a manner analogous to that of functions of a complex variable. The power of the latter, key to much of the power of mathematical analysis, depends largely on restricting attention to a particular class of functions of a complex variable, namely the holomorphic functions. Possibly the reason for Hamilton's neglect of this approach was because the most obvious generalisation of the complex holomorphic condition is the existence, for function $f$, of a function $g$ such that $df = gdz$ where $z$ is the independent complex variable. Interpreting these quantities instead as quaternions, this condition becomes a constraint so severe that only functions linear in all four coordinates of $z$ can meet it (and only a subset of those). Such a class of functions is so narrow as to be of little interest.

Apparently it was Fueter (1935,1936) who first realised that the generalisation of the concept of a holomorphic function into the quaternionic domain was more productively realised by generalising the Cauchy-Riemann conditions and by noting that $(\partial_x + i\partial_y)f = 0$ for the complex domain is extended into the quaternionic domain as $(\partial_w + I\partial_x + J\partial_y + K\partial_z)f = 0$, thus defining "regular" quaternionic functions.

Fueter derived an integral theorem respected by such regular functions and which can be regarded as analogous to Cauchy's Theorem for holomorphic complex functions. The natural expression of such integrals, essentially geometrical integrals, is in terms of exterior forms so that all such integral theorems, including Cauchy's Theorem, are derived from the Stokes-Cartan theorem, $\oint_{\partial C} F = \int_C dF$. In the quaternion domain, the manifold $C$ is 4-dimensional, which implies that $F$ must be a 3-form. This is why it will not do to impose $df = gdz$ with $g$ a function (a 0-form) as it implies that the relevant $F = df$ is a 1-form. Instead, the property of regularity emerges as equivalent to a definition of the "derivative" as an identity between two 3-forms (see below).

Useful review articles have been published by Deavours (1973), Sudbery (1979), and Fokas and Pinotsis (2007). These discuss the derivation of Fueter's Theorem, as well as the conditions under which the Taylor or Laurent series exist, together with a corresponding residue theorem and a generalisation of Liouville's Theorem.

The purpose of the present paper is to derive a new class of integral theorems for regular quaternionic functions. At first sight these may appear to be but a small step away from Fueter's Theorem, and hence not fundamentally new. However, on closer inspection it is not clear (to this author, at least) how one theorem can be derived directly from the other. Moreover, the key property of both theorems is their independence of integration surface, and this boundary independence is shown to arise in different ways for the two theorems.

The opportunity is taken to discuss the extension into the biquaternion domain of both this new integral theorem and Fueter's Theorem. The relevance of this in physics is that the Lorentz transformation of Minkowski spacetime is elegantly expressed in terms of Hermitian biquaternions, Lambek (2013). Consequently, functions which obey the corresponding extension of regularity, namely $(\partial_w + i[I\partial_x + J\partial_y + K\partial_z])f = 0$, obey the wave equation. In the biquaternion case, integration necessary takes place on the Argand plane. The integration contour will require careful definition with respect to which poles of the integrand are included and which excluded from the closed contour.

## Quaternion and Biquaternion Notation

The general quaternion is composed of four real components, or coordinates, $w, x, y, z$,

$$q = w + xI + yJ + zK \qquad (1)$$

where the basis quaternions obey $I^2 = J^2 = K^2 = IJK = -1$, and hence are non-commutative, e.g, $IJ = -JI = K$, etc. The general biquaternion can also be written as (1) but the four coordinates $w, x, y, z$ may now be complex valued. The quaternion-conjugate is defined as,

$$q^\# = w - xI - yJ - zK \qquad (2)$$

which contrasts with the usual complex conjugate, which is,

$$q^* = w^* + x^*I + y^*J + z^*K \qquad (3)$$

A quaternion may be written $q = v_t + \bar{v}$, where $v_t$ is a real number (or a complex number in the biquaternion case) and $\bar{v} = v_x I + v_y J + v_z K$ in which the basis quaternions play the part of a coordinate basis of unit vectors. Hence, the part of a quaternion dependent upon $I, J, K$ is referred to as the vector part whilst the rest is the scalar part. The product of two quaternions is,

$$ab = (a_w b_w - \bar{a} \cdot \bar{b}) + (a_w \bar{b} + b_w \bar{a} + \bar{a} \times \bar{b}) \tag{4}$$

and hence produces the usual scalar (dot) product and the vector (cross) product of vectors automatically as a consequence of the quaternion algebra, as well as the multiplication of a vector by a scalar. (4) gives $(ab)^\# = b^\# a^\#$. Quaternions have the natural squared-norm,

$$\mathcal{N}(a) = aa^\# = a^\# a = a_w^2 + |\bar{a}|^2 = a_w^2 + a_x^2 + a_y^2 + a_z^2 \tag{5}$$

All quaternions except zero therefore have a positive-definite squared-norm, whilst this is not the case for biquaternions. Quaternions therefore form a division ring because all non-zero quaternions have an inverse, which is given explicitly by $a^{-1} = a^\# / \mathcal{N}(a)$. This is precluded for biquaternions because $\mathcal{N}(a)$ may be zero, and hence biquaternions do not form a division ring, some biquaternions having no inverse, e.g., $1 + iI$ does not (where $i$ has its usual interpretation in complex numbers as the commuting form of $\sqrt{-1}$).

In terms of a unit vector, $\hat{n} = In_x + Jn_y + Kn_z$ where $n_x^2 + n_y^2 + n_z^2 = 1$, an arbitrary quaternion can be written in "polar" form,

$$a = re^{\theta \hat{n}} = r(\cos(\theta) + \hat{n} \sin(\theta)) \tag{6}$$

and this prescription is unique within the ranges $r \geq 0$, $0 \leq \theta < \pi$ if $\hat{n}$ is allowed to roam over all eight octants. The original reason for Hamilton's interest in quaternions was that rotations in Euclidean $\mathbb{R}^3$ can be written $p \to \tilde{p} = qpq^{-1}$ for an arbitrary quaternion, $p$, in terms of the unit quaternion $q = e^{\theta \hat{n}}$, and hence $q^{-1} = e^{-\theta \hat{n}}$, where $\hat{n}$ is the axis of rotation and $\theta$ is the angle of rotation.

A particular class of biquaternions are the Hermitian biquaternions, defined by $q^* = q^\#$, or, equivalently, $q^{\#*} = q$. Hermitian biquaternions can be written in terms of four real coordinates, $q_\mu$, as $q = q_t + i\bar{q}$, where $\bar{q} = q_x I + q_y J + q_z K$. The squared-norm of an Hermitian biquaternion is,

$$\mathcal{N}(q) = qq^\# = q^\# q = q_t^2 - |\bar{q}|^2 = q_t^2 - (q_x^2 + q_y^2 + q_z^2) \tag{7}$$

Hence, unlike quaternions which have a natural Euclidean norm, (5), the Hermitian biquaternions have a natural Minkowskian norm, (7). Note that in (7) we have changed the notation for the fourth coordinate from $w$ to $t$ in recognition that the physical interpretation as time is now possible. Consequently, events in Minkowskian spacetime may be represented as Hermitian biquaternions. It follows that the transformation of an arbitrary Hermitian biquaternion, $p$, defined by,

$$p \to p' = qpq^{*\#} \tag{8}$$

where $q$ is any biquaternion with unit norm (i.e., with $q^\# q = qq^\# = 1$) is a Lorentz transformation. If we write $q = u + iv$, where $u$ and $v$ are (real) quaternions, then $q$ is a rotation if it is a quaternion (i.e., real, hence $v = 0$) whereas $q$ is a boost if it is Hermitian, i.e., if $u$ is purely scalar and $v$ is purely vector, i.e., $\bar{u} = 0, v_t = 0$.

Hermitian biquaternions have a natural expression as exponentials. Any spacetime event at a timelike 4-vector displacement from the origin can be represented as $p = t + i\bar{r}$ and written as,

$$p = re^{i\theta \hat{n}} = r(\cosh(\theta) + i\hat{n} \sinh(\theta)) \tag{9}$$

Similarly, any spacetime event at a spacelike 4-vector displacement from the origin can be represented as

$$p = ire^{\left(\frac{\pi}{2}+i\theta\right)\hat{n}} = r\bigl(sinh(\theta) + i\hat{n}cosh(\theta)\bigr) \qquad (10)$$

where the real (quaternionic) unit vector, $\hat{n}$, is allowed to roam over all eight octants. All non-null spacetime events can be uniquely represented in the form (9) or (10) within the range $0 \leq \theta < \infty$ with $r$ allowed to take any value on the real line, positive or negative. Events at a null (light-like) displacement from the origin are obtained in the limit $\theta \to \infty$.

## Regular Functions

Until further notice we confine attention to quaternions, i.e., with real coordinates.

Complex holomorphic functions are defined by the existence of a derivative, so that the complex valued function $f$ is said to be holomorphic in the neighbourhood of a point $z$ in the Argand plane if there exists $g$ such that,

$$df = gdz \qquad (11)$$

This has the implication that the quotient $df/dz$ is the same, namely $g$, irrespective of the direction from which the limiting point is approached. Writing $f = f_x + if_y$ where $f_x$ and $f_y$ are real, this identifies holomorphic functions as those respecting the Cauchy-Riemann equations,

$$f_{x,x} = f_{y,y} \quad \text{and} \quad f_{x,y} = -f_{y,x} \qquad (12)$$

and hence that holomorphic functions obey Laplace's equation in $\mathbb{R}^2$,

$$\nabla_2^2 f \equiv \bigl(\partial_x^2 + \partial_y^2\bigr)f = 0 \qquad (13)$$

The failure of the nineteenth century mathematicians to develop the analysis of quaternionic functions in an analogous manner might be due to the following fact. Suppose we define a class of quaternion valued functions in a manner directly analogous to (11), i.e., that a unique derivative, $g$, exists with respect to the change in quaternionic coordinates, $dq$, such that $df = gdq$ irrespective of the direction in $\mathbb{R}^4$ from which the limiting point in $\mathbb{R}^4$ is approached. This is readily seen to be a constraint so severe that only functions linear in all the coordinates of $q$ can meet it (and only a subset of those). Such a class of functions is so narrow as to be of little interest.

Apparently, it was Fueter (1935,1936) who first realised that the generalisation of the concept of a holomorphic function into the quaternionic domain was more productively realised by generalising the Cauchy-Riemann conditions, (12), and by noting that these can be written as a single, complex valued, equation as $\bigl(\partial_x + i\partial_y\bigr)f = 0$. Hence, the quaternionic differential operator, $D$, is defined by analogy as,

$$D = \partial_w + I\partial_x + J\partial_y + K\partial_z \equiv \partial_w + \overline{\nabla} \qquad (14)$$

where $\overline{\nabla}$ is understood to be the usual vector operator but in quaternionic form. Hence a quaternion valued function of the four real variables $w, x, y, z$ is said to be "regular" in some domain of $\mathbb{R}^4$ if throughout that region it obeys,

$$Df = 0 \qquad (15)$$

The term "regular" is used rather than "holomorphic" to distinguish condition (15) from condition $df = gdq$. Unlike the complex case, the latter condition does not follow from the former. We shall shortly present the correct quaternion analogue of $df = gdq$.

The quaternion conjugate operator $D^\# \equiv \partial_t - \overline{\nabla}$ is such that the Laplacian in $\mathbb{R}^4$ can be written $\nabla_4^2 = D^\# D = DD^\#$. Hence it immediately follows from the definition, (15), that regular quaternion functions obey Laplace's equation in $\mathbb{R}^4$, $\nabla_4^2 f = 0$. This generalises the naturally harmonic nature of holomorphic complex functions.

We can also define conjugate-regular functions by $D^\# g = 0$, which also obey $\nabla_4^2 g = 0$. Note that $D^\# g$ does not in general equal either $(Dg)^\#$ or $(Dg^\#)^\#$.

Regular functions will also be called "left regular" because we have,

$$Df = \partial_t f + I\partial_x f + J\partial_y f + K\partial_z f = 0 \tag{16}$$

i.e., the basis quaternions stand on the left of the function. We can also define for later purposes the right-regular functions which obey,

$$\widetilde{D}f \equiv \partial_t f + \partial_x f I + \partial_y f J + \partial_z f K = 0 \tag{17}$$

Note that an alternative notation is $\widetilde{D}f \equiv f\overleftarrow{D}$, where the arrow denotes that the derivatives operate on the function to the left, whereas, for the purposes of the quaternion content, the operator stands on the right. We shall see that both left and right regular functions play a key role in the most cogent explanation of the origin of Fueter's Theorem.

Similarly, conjugate-right-regular functions obey,

$$\widetilde{D}^\# f \equiv \partial_t f - \partial_x f I - \partial_y f J - \partial_z f K = 0 \tag{18}$$

So we have $(Dg)^\# \equiv \widetilde{D}^\# g^\#$ and $(\widetilde{D}g)^\# \equiv D^\# g^\#$. Right-regular functions and conjugate right-regular functions also obey Laplace's equation, $\nabla_4^2 g = 0$.

All the above differential operators, and the definition of "regular" functions, are appropriate for use with quaternion-valued functions, i.e., with real coordinates. Their extension to biquaternions, and hence to the Minkowski metric, will be addressed later.

**Quaternion Geometrical Integrals**

It is instructive to recall how the derivation of the holomorphic integral theorem for complex functions arises from the Stokes-Cartan Theorem. The latter states that,

$$\oint_{\partial C} F = \int_C dF \tag{19}$$

where $F$ is a smooth $(n-1)$-form and $d$ is to be understood as the exterior derivative, so that $dF$ is a $n$-form. In (19) $C$ is a simply connected orientable $n$-dimensional manifold and $\partial C$ is its $(n-1)$ dimensional boundary. Hence, for $n = 2$ and $C$ a planar Euclidean manifold, (19) becomes Green's Theorem, whilst if $n = 2$ and the 2-surface $C$ is embedded within a three-dimensional Euclidean manifold, (19) yields the original Stokes' Theorem. Using $n = 3$ when $C$ is a region of a three-dimensional Euclidean manifold, (19) produces Gauss's Theorem (the divergence theorem in 3D).

The holomorphic integral theorem is what Green's Theorem implies for a complex function, $f$, for which (11) holds. This is sufficient to imply the Cauchy-Riemann conditions, (12), which can also be written $D_c f = 0$ where $D_c = (\partial_x + i\partial_y)$. The 1-form to be inserted into the LHS of (19) is $F = fdz = fdx + ifdy$. Hence on the RHS of (19),

$$d(fdz) = df \wedge dx + idf \wedge dy = (\partial_y f)dy \wedge dx + i(\partial_x f)dx \wedge dy = i(D_c f)dx \wedge dy = 0 \quad (20)$$

Hence, (19) becomes the holomorphic integral theorem $\oint_{\partial C} fdz = 0$ and Cauchy's Theorem immediately follows. In order to generalise this approach to quaternions we shall need an integral which encloses a region in $\mathbb{R}^4$, i.e., an integral over a 3-surface. Hence, if (19) is to be deployed, the $F$ on the LHS must be a 3-form and the resulting 4-form $dF$ on the RHS can only be proportional to the natural volume element, $V_q = dw \wedge dx \wedge dy \wedge dz$. What, then, can the 3-form $F$ be if it is to depend only upon the quaternion function, $f$, and the quaternion manifold itself?

Hence we see why the sensible analogue of the "derivative" of a complex function, namely $g$ in $df = gdz$ (a 1-form), cannot for quaternions be $df = gdq$ because this is also a 1-form whereas a 4-form is needed, namely a function (a 0-form) times $V_q = dw \wedge dx \wedge dy \wedge dz$. The solution to this conundrum is that the sensible (or productive) definition of the "derivative" of a quaternionic function is the $g$ in,

$$\frac{1}{2} d(dq \wedge dq f) = S_q g \quad (21)$$

where $S_q$ is the natural 3-surface element in quaternion space, namely,

$$S_q = dx \wedge dy \wedge dz - Idw \wedge dy \wedge dz - Jdw \wedge dz \wedge dx - Kdw \wedge dx \wedge dy \quad (22)$$

and $q$ is the position quaternion, (1), so that the 1-form $dq$ is the identity mapping. The four terms in $S_q$ can be interpreted as the projections of the 3-surface onto the axes $w, x, y, z$ respectively. The wedge product $dq \wedge dq$ is a little startling at first because it would be identically zero for commuting coordinates. This is not so in the quaternionic case because quaternions do not commute. Instead we have,

$$\frac{1}{2} dq \wedge dq = \frac{1}{2}(dw + Idx + Jdy + Kdz) \wedge (dw + Idx + Jdy + Kdz)$$
$$= Idy \wedge dz + Jdz \wedge dx + Kdx \wedge dy \quad (23)$$

which can be interpreted as the natural 2-surface element within the 3-dimensional $x, y, z$ spatial subspace of quaternion space. The imposition of (21) as a condition on $f$ is equivalent to the definition of a regular function via (14,15), i.e., being regular ensures that a "derivative", $g$, satisfying (21) exists, and vice-versa. This is shown by explicit evaluation of the LHS of (21) as follows,

$$\frac{1}{2} d(dq \wedge dq f) = d(Idy \wedge dz f + Jdz \wedge dx f + Kdx \wedge dy f)$$
$$= Idy \wedge dz \wedge (f_{,w} dw + f_{,x} dx) + Jdz \wedge dx \wedge (f_{,w} dw + f_{,y} dy) + Kdx \wedge dy \wedge (f_{,w} dw + f_{,z} dz)$$
$$= (Idw \wedge dy \wedge dz + Jdw \wedge dz \wedge dx + Kdw \wedge dx \wedge dy) f_{,w} + dx \wedge dy \wedge dz (If_{,x} + Jf_{,y} + Kf_{,z})$$
$$= (Idw \wedge dy \wedge dz + Jdw \wedge dz \wedge dx + Kdw \wedge dx \wedge dy) f_{,w} + dx \wedge dy \wedge dz \overline{\nabla} f \quad (24)$$

If we assume that $f$ is regular, so that (15) requires $Df = (\partial_w + \overline{\nabla})f = 0$, then the expression (24) becomes,

$$[(Idw \wedge dy \wedge dz + Jdw \wedge dz \wedge dx + Kdw \wedge dx \wedge dy) - dx \wedge dy \wedge dz]f_{,w} = -S_q f_{,w} \qquad (25)$$

i.e., (21) is established with $g = -f_{,w} = \overline{\nabla}f$, the "derivative" of $f$. Conversely, if (21) is taken as the defining condition, then requiring that (24) is an identity with $S_q g$ requires the regularity condition, (15).

In the same way it is readily checked that a conjugate-regular function obeying $D^\# f = (\partial_w - \overline{\nabla})f = 0$ is equivalent to $\frac{1}{2}d(dq \wedge dqf) = S_q^\# g$ and hence to the existence of a "derivative", $g = f_{,w} = \overline{\nabla}f$. Regular and conjugate-regular functions are thus interchangeable by reversal of the Euclidean "time" coordinate, $w$.

That the 2-surface element $\frac{1}{2}dq \wedge dq$ stands to the left of $f$ in (21) is crucial to the above results. Reversing the order is readily seen to provide a condition equivalent to the function being right-regular, (17), i.e.,

$$\tfrac{1}{2}d(fdq \wedge dq) = gS_q \Leftrightarrow \widetilde{D}f = 0 \qquad (26)$$

And similarly conjugate-right-regular functions, (18), have $\frac{1}{2}d(fdq \wedge dq) = gS_q^\#$.

**The Quaternionic Analogue of the Holomorphic Integral Theorem**

In the complex domain the holomorphic integral theorem is what results from (19,20), which explicitly is,

$$\oint_{\partial C} fdz = 0 \qquad (27)$$

for any closed contour, $\partial C$, in the Argand plane and for any function which is holomorphic everywhere within, and on, $\partial C$. As we now know that the natural quaternionic 3-form (or 3-surface element) is $S_q$, (22), the quaternion analogue of (27) is expected to be,

$$\oint_{\partial V} S_q f = 0 \qquad (28)$$

where $\partial V$ is the boundary of any simply connected region, $V$, in quaternion $\mathbb{R}^4$ space and $f$ is any function which is regular everywhere within, and on, $\partial V$. The order of the two terms in (28) is important because they do not commute. (Reversing the order will be considered below). To establish (28), note that the Stokes-Cartan Theorem, (19), tells us that,

$$\oint_{\partial V} S_q f = \int_V d(S_q f) \qquad (29)$$

Hence (28) is immediately established because $d(S_q f) = 0$ for a regular function, which can be seen as follows, using, from (22), that $dS_q = 0$,

$$d(S_q f) = S_q df$$
$$= (dx \wedge dy \wedge dz - Idw \wedge dy \wedge dz - Jdw \wedge dz \wedge dx - Kdw \wedge dx \wedge dy) \wedge (f_{,w}dw + f_{,x}dx + f_{,y}dy + f_{,z}dz)$$
$$= -dw \wedge dx \wedge dy \wedge dz(f_{,w} + If_{,x} + Jf_{,y} + Kf_{,z}) = -V_q Df \qquad (30)$$

which is zero for a regular function because $Df = 0$. QED.

Similarly, for conjugate-regular functions, $\oint_{\partial V} S_q^\# f = 0$ because $\oint_{\partial V} S_q^\# f = \int_V d(S_q^\# f)$ and $d(S_q^\# f) = -V_q D^\# f = 0$. And for right-regular functions, $\oint_{\partial V} fS_q = 0$ because $d(fS_q) =$

$(\widetilde{D}f)V_q = 0$, and so, for conjugate-right-regular functions, $\oint_{\partial V} fS_q^{\#} = 0$. Note the reversal of the order of the two quaternionic terms for the right-regular cases.

## Existence and Quaternionic Analysis

We have not been explicit about the conditions under which all the above quantities exist. It is not the purpose of this brief note to address this in detail, nor the relevance of these conditions to broader aspects of quaternionic analysis. However a few remarks are appropriate.

It is sufficient to ensure that the Stokes-Cartan Theorem holds, and hence the existence of the integrals derived herein, to assume that all partial derivatives exist and are continuous. However, these results hold in more general circumstances. Sudbery (1979), quoting Schuler (1937) as the source, claims that (28), and by implication the generalisations of Cauchy's Theorem in the following sections, only requires the existence of the partial derivatives, though perhaps restricted to rectifiable boundaries. Sudbery goes on to claim that,

*"From this it follows, as in complex analysis, that if f is regular in an open set U then it has a power series expansion about each point of U. Thus, pointwise differentiability, together with the four real conditions, (15), on the sixteen partial derivatives of f, is sufficient to ensure analyticity."*

Hence, a vista opens in which quaternionic analysis may be as rich as complex analysis, supporting Taylor or Laurent series, and an associated residue theorem, as well as a generalisation of Liouville's Theorem. The latter establishes that any regular quaternionic function which is bounded across the whole of quaternion $\mathbb{R}^4$ space is a constant. For further details on these matters see Deavours (1973) or Fokas and Pinotsis (2007).

## Fueter's Theorem

We are now interested in obtaining quaternionic generalisations of Cauchy's Theorem for complex holomorphic functions. The latter is,

$$\oint_{\partial C} \frac{f}{z-z_0} dz = 2\pi i f(z_0) \tag{31}$$

where $\partial C$ is any closed contour which fully contains the point $z_0$ and $f$ is any function which is holomorphic everywhere inside, and on, $\partial C$. The contour independence of the integral, providing it continues to fully enclose $z_0$, follows immediately from (27) because the integrand $\frac{f}{z-z_0}$ is holomorphic everywhere other than $z_0$. The value of the integral then follows by explicit evaluation on a limiting circular contour shrinking onto point $z_0$.

To extend to the case of regular quaternion functions, by analogy this means finding an appropriate auxiliary function, $H$, which must be singular at a single isolated point, $q_0$, to factor into the integrand of $\oint_{\partial V} S_q f$ so that,

(i) on a limiting surface which shrinks onto $q_0$ the integral evaluates to a finite, non-zero numerical constant times $f(q_0)$, and,

(ii) the integral is independent of the boundary, $\partial V$, provided that point $q_0$ is fully contained.

The latter requirement is equivalent to the modified integral being zero for any closed surface which fully excludes $q_0$.

The distinction between Fueter's Theorem, to be presented in this section, and the theorem which is the subject of this paper, to be presented in the following section, lies in the different approaches to how condition (ii) is achieved. Presentations of Fueter's Theorem often do not spell out the reason for the integration being surface independent. Fueter's Theorem fulfils the two conditions using the following integral, for a left-regular function, $f$,

$$\oint_{\partial V} \frac{(q-q_0)^{-1}}{\mathcal{N}(q-q_0)} S_q f(q) = 2\pi^2 f(q_0) \qquad (32)$$

(We have taken liberties with the notation here, and throughout, in that a regular function is not really just a function of $q$ so $f(q)$ is actually a short-hand for $f(w, x, y, z)$).

The value of the constant can be found by considering $\partial V$ to be a 3-sphere of radius $r$ centred on $q_0$, with the limit $r \to 0$ in mind. The auxiliary function $\frac{(q-q_0)^{-1}}{\mathcal{N}(q-q_0)}$ is of order $\frac{1}{r^3}$ which cancels with the $r$ dependence of the surface element 3-form, $S_q$, so the resulting integral is finite and non-zero. Detailed evaluation shows that it reduces simply to the surface area of a 3-sphere of unit radius, i.e., $2\pi^2$.

However, the interesting feature of (32) is how its surface independence arises. This is due to the following theorem: if $g$ is a right-regular function and $f$ is a left-regular function, then the following integral is identically zero,

$$\oint_{\partial V} g(q) S_q f(q) = 0 \qquad (33)$$

independent of the closed boundary, $\partial V$. This is established by showing that the exterior derivative of the integrand is zero, $d(gS_q f) = 0$, and once again applying the Stokes-Cartan Theorem, (19). Note that the order of the terms in (33) is crucial; the 3-form surface element, $S_q$, must be sandwiched between the two functions, as will be clear from the proof - which is simple, as follows.

$$d(gS_q f) = (dg) \wedge S_q f + gS_q \wedge (df)$$

$$= (g_{,w} dw \wedge dx \wedge dy \wedge dz - g_{,x} I dx \wedge dw \wedge dy \wedge dz - g_{,y} J dy \wedge dw \wedge dz \wedge dx - g_{,z} K dz \wedge dw \wedge dx \wedge dy) f +$$
$$g(dw \wedge dx \wedge dy \wedge dz f_{,w} - I dx \wedge dw \wedge dy \wedge dz f_{,x} - J dy \wedge dw \wedge dz \wedge dx f_{,y} - K dz \wedge dw \wedge dx \wedge dy f_{,z})$$

$$= (g_{,w} + g_{,x} I + g_{,y} J + g_{,z} K) V_k f + gV_k (f_{,w} + I f_{,x} + J f_{,y} + K f_{,z})$$

$$= (\tilde{D} g) V_k f + gV_k D f = 0 \qquad (34)$$

The last step follows from the assumption that $g$ is right-regular, $\tilde{D} g = 0$, and that $f$ is left-regular, $Df = 0$. Note that the basis quaternions, $I, J, K$, must be kept sandwiched between the $g$ and $f$ terms throughout in order for (34) to hold, hence the importance of the order of the three terms in (33).

It remains to be shown that the function,

$$H(q - q_0) = \frac{(q-q_0)^{-1}}{\mathcal{N}(q-q_0)} = \frac{(q-q_0)^{\#}}{|q-q_0|^4} \qquad (35)$$

which appears in (33) is left-regular. In actual fact this function is unusual in being both right-regular and left-regular. The simplest demonstration of regularity is to note that this function can be written,

$$H(q) = \frac{q^\#}{|q|^4} = -e^{-w\bar{\nabla}}\left(\frac{\bar{r}}{r^4}\right) \tag{36}$$

where $q$ is given by (1), i.e., $q = w + \bar{r}$, so that $r^2 = x^2 + y^2 + z^2$. (36) can be shown by first showing $\bar{\nabla}\bar{r} = -3$, and for integral $n$ that $\bar{\nabla}\left(\frac{1}{r^n}\right) = -n\frac{\bar{r}}{r^{n+2}}$, whilst for $n = 1, 3, 5, \ldots$ $\bar{\nabla}\left(\frac{\bar{r}}{r^{n+3}}\right) = \frac{n}{r^{n+3}}$. Expanding the exponential in (36) and summing the resulting series provides the result. That $H$ is left-regular now follows immediately because,

$$DH = -(\partial_w + \bar{\nabla})e^{-w\bar{\nabla}}\left(\frac{\bar{r}}{r^4}\right) = -(-\bar{\nabla} + \bar{\nabla})e^{-w\bar{\nabla}}\left(\frac{\bar{r}}{r^4}\right) = 0 \tag{37}$$

That $H$ is also right-regular is seen by explicit evaluation of the derivatives using $H(q) = \frac{q^\#}{|q|^4}$ which gives,

$$\partial_w H = -H\overleftarrow{\nabla} = \frac{r^2 + 4w\bar{r} - 3w^2}{|q|^6} \quad \text{hence} \quad \partial_w H + H\overleftarrow{\nabla} = \widetilde{D}H = 0 \tag{38}$$

This completes the proof of Fueter's Theorem, (32).

It immediately follows that, for a right-regular function, $f$,

$$\oint_{\partial V} f(q) S_q \frac{(q-q_0)^{-1}}{N(q-q_0)} = 2\pi^2 f(q_0) \tag{39}$$

which exploits the fact that $H$ is also left-regular. For conjugate-regular functions, $f$, we have,

$$\oint_{\partial V} \frac{q-q_0}{|q-q_0|^4} S_q^+ f(q) = 2\pi^2 f(q_0) \tag{40}$$

whilst for conjugate-right-regular functions, $f$,

$$\oint_{\partial V} f(q) S_q^+ \frac{q-q_0}{|q-q_0|^4} = 2\pi^2 f(q_0) \tag{41}$$

**An Alternative Cauchy-like Integral Theorem**

The use of theorem (34) to establish a Cauchy-like integral formula is not the most obvious approach to satisfying condition (ii). More obvious would be to consider the product of the auxiliary function, $H$, with the regular function, $f$, and to seek an integral of the form $\oint_{\partial V} S_q H(q - q_0) f(q)$. This would immediately ensure boundary-surface independence if the product $H(q - q_0) f(q)$ were regular everywhere except at the singular point, $q_0$. In the case of complex holomorphic functions, this approach works because of the chain rule of differentiation. That is,

$$(\partial_x + i\partial_y)(fg) = [(\partial_x + i\partial_y)f]g + f(\partial_x + i\partial_y)g$$

so the product of two holomorphic functions is holomorphic. For regular quaternion functions this fails because the chain rule fails due to quaternions not commuting. Thus we cannot in general write $D(fg)$ as $(Df)g + fDg$ because the quaternionic operator cannot be moved to the right of the quaternionic function $f$. No doubt this is what discouraged Fueter and co-workers from employing this approach.

However, there are classes of regular function for which the product is regular. The most obvious is when one of the functions is purely scalar, and hence the functions commute. However, this is quite useless as a purely scalar function which is also regular must be a constant.

However, there is another class of regular function which proves more useful. It is convenient to note, following Hamilton and Graves, that any regular function, $g(w, x, y, z)$, can be written in terms of a "generating function", $G(x, y, z)$, as follows,

$$g = e^{-w\bar{\nabla}} G \tag{42}$$

That $Dg = 0$ follows immediately for any differentiable $G(x, y, z)$, which is, of course, simply $g(0, x, y, z)$. Both $g$ and $G$ are quaternionic in general.

We can now prove this theorem: If $g$ and $f$ are both regular and $g$ has a scalar generating function, then $gf$ is regular. Writing $f = e^{-w\bar{\nabla}} F$ together with $g = e^{-w\bar{\nabla}} G$, where $G$ is scalar, we have,

$$D(gf) = (\partial_w + \bar{\nabla})\left[\left(e^{-w\bar{\nabla}} G(x,y,z)\right)\left(e^{-w\bar{\nabla}} F(x,y,z)\right)\right] \tag{43}$$

We note that $\partial_w$ commutes with $G$ and $F$ but not with the factors $e^{-w\bar{\nabla}}$, whereas $\bar{\nabla}$ commutes with the factors $e^{-w\bar{\nabla}}$ but not with $G$ or $F$ when these are quaternionic Carrying out the $w$-derivatives gives,

$$D(gf) = \left(-\bar{\nabla} e^{-w\bar{\nabla}} G(x,y,z)\right)\left(e^{-w\bar{\nabla}} F(x,y,z)\right) + \left(e^{-w\bar{\nabla}} G(x,y,z)\right)\left(-\bar{\nabla} e^{-w\bar{\nabla}} F(x,y,z)\right) \\ + \bar{\nabla}\left[\left(e^{-w\bar{\nabla}} G(x,y,z)\right)\left(e^{-w\bar{\nabla}} F(x,y,z)\right)\right] \tag{44}$$

We now impose the condition that $G$ is scalar. In the last term of (44) the quaternionic aspects of $\bar{\nabla}$ now commute with $G$ which allows the simple chain rule to hold, $\bar{\nabla}(GF) = (\bar{\nabla}G)F + G(\bar{\nabla}F)$, with the result that (44) is zero. QED.

This provides another route to a boundary-independent integral, using $\oint_{\partial V} S_q g(q - q_0) f(q)$ where $g$ has a scalar generating function, $G$. But what is a suitable $G$ to satisfy conditions (i) and (ii)?

We see that the function $H(q) = \frac{q^\#}{|q|^4}$ deployed in the Fueter integral will not do because its generating function is $-\frac{\bar{r}}{r^4}$, (36), which is not scalar.

Considering the integral to be evaluated on a 3-sphere of radius $\rho$ centred on $q_0$ shows that any candidate function must be of the order of $1/\rho^3$ as $\rho \to 0$ in order to produce an integral which is finite and non-zero. The most obvious candidate, $g = (q - q_0)^{-3}$, is not regular. The regular function formed from it by using a generating function $G = (\bar{r} - \bar{r}_0)^{-3}$ will not do either as this is not scalar (and just reproduces Fueter's $-H(q) = -\frac{q^\#}{|q|^4}$).

The next most obvious candidate is to use the generating function $G = |\bar{r} - \bar{r}_0|^{-3}$ which is scalar and hence will produce a boundary-independent integral. However, despite having a singularity of the required order at $q_0$, the integral is identically zero even when the integration surface contains $q_0$, because the angular integrations make it zero. This may be checked directly noting that the auxiliary function is then,

$$g = e^{-w\bar{\nabla}}\left(\frac{1}{r^3}\right) = \frac{r^2 - w^2}{r(r^2 + w^2)^2} + \frac{w(3r^2 + w^2)}{r^3(r^2 + w^2)^2}\bar{r} \tag{45}$$

In order to prevent the angular integrations being zero, an angular dependence needs to be introduced into the generating function whilst ensuring it remains scalar and has a singularity of the order of $1/\rho^3$ as $\rho \to 0$. The simplest generating function meeting these requirements is $\frac{x}{r^4}$, or the equivalent with $x$ replaced by $y$ or $z$. Hence we explore the option,

$$g = e^{-w\bar{\nabla}}\left(\frac{x}{r^4}\right) \tag{46}$$

Explicit evaluation is surprising as this function turns out to be transcendental, (47)

$$g = e^{-w\bar{\nabla}}\left(\frac{x}{r^4}\right) = \frac{x}{(w^2+r^2)^2} - I\frac{w}{2r^2}\left\{\frac{1}{w^2+r^2} + \frac{1}{rw}tan^{-1}\left(\frac{w}{r}\right)\right\} + 4xw\frac{\bar{r}}{r^4}\left\{\frac{5r^2+3w^2}{8(w^2+r^2)^2} + \frac{3}{8rw}tan^{-1}\left(\frac{w}{r}\right)\right\}$$

However, because the generating function is scalar, the integral $\oint_{\partial V} S_q g(q-q_0) f(q)$ is ensured to be boundary independent for any regular function $f$ where $\partial V$ contains $q_0$, and hence the integral is a constant times $f(q_0)$. The constant can be evaluated on any convenient boundary. Choosing either a 3-sphere centred on $q_0$ or the limit of an infinitely long prismatic surface oriented parallel to the $w$-axis, and with a cross-section which is a 2-sphere, shows the constant to be $\frac{2\pi^2}{3}I$.

Hence, finally, we have our alternative Cauchy-like integral theorem for any regular quaternion function, $f$,

$$\oint_{\partial V} S_q \left[e^{-(w-w_0)\bar{\nabla}}\left(\frac{x-x_0}{|\bar{r}-\bar{r}_0|^4}\right)\right] f(q) = \frac{2\pi^2}{3} If(q_0) \tag{48}$$

Unlike the Fueter integral, the constant in (48) is factored by the basis quaternion, $I$. In this respect it is more like the Cauchy integral whose constant also involves a square root of -1, i.e., $2\pi i$. By symmetry we also have,

$$\oint_{\partial V} S_q \left[e^{-(w-w_0)\bar{\nabla}}\left(\frac{y-y_0}{|\bar{r}-\bar{r}_0|^4}\right)\right] f(q) = \frac{2\pi^2}{3} Jf(q_0) \tag{49}$$

$$\oint_{\partial V} S_q \left[e^{-(w-w_0)\bar{\nabla}}\left(\frac{z-z_0}{|\bar{r}-\bar{r}_0|^4}\right)\right] f(q) = \frac{2\pi^2}{3} Kf(q_0) \tag{50}$$

At first sight it seems that the Fueter integral should be simply derivable from (48-50), by multiplying these equations by $I, J$ or $K$ respectively and summing the result. This looks promising because,

$$-e^{-(w-w_0)\bar{\nabla}}\left(\frac{(x-x_0)I+(y-y_0)J+(z-z_0)K}{|\bar{r}-\bar{r}_0|^4}\right)$$

is just the Fueter function $H(q-q_0)$, (36). The problem with carrying out this programme is that the basis quaternions, $I, J$ or $K$, cannot be pushed through the $S_q$ or the $e^{-(w-w_0)\bar{\nabla}}$ or the $f(q)$ terms in the integrals (48-50) due to lack of commutation. Moreover, to reproduce the Fueter integral, (32), the function $H(q-q_0)$ must stand to the left of $S_q$. Perhaps the reader can see a way of deriving Fueter, (32), from (48-50), but it is not obvious to the present author.

However, more importantly, it seems even less likely that the reverse derivation can be done. Bearing in mind that the boundary independence of the Fueter integral, (32), and that of (48-50) rely on quite different properties, it seems that (48-50) are genuinely distinct theorems. In some respects (48-50) are more akin to Cauchy than (32) in that the product $\left[e^{-(w-w_0)\bar{\nabla}}\left(\frac{x-x_0}{|\bar{r}-\bar{r}_0|^4}\right)\right] f(q)$ is regular except at $q_0$.

We derive in the same way the following,

For conjugate regular $f$,

$$\oint_{\partial V} S_q^{\#} \left[ e^{(w-w_0)\bar{\nabla}} \left( \frac{x-x_0}{|\bar{r}-\bar{r}_0|^4} \right) \right] f(q) = \frac{2\pi^2}{3} I f(q_0) \tag{51}$$

For right-regular $f$,

$$\oint_{\partial V} f(q) \left[ \left( \frac{x-x_0}{|\bar{r}-\bar{r}_0|^4} \right) e^{-(w-w_0)\bar{\nabla}} \right] S_q = \frac{2\pi^2}{3} f(q_0) I \tag{52}$$

For conjugate right-regular $f$,

$$\oint_{\partial V} f(q) \left[ \left( \frac{x-x_0}{|\bar{r}-\bar{r}_0|^4} \right) e^{(w-w_0)\bar{\nabla}} \right] S_q^{\#} = \frac{2\pi^2}{3} f(q_0) I \tag{53}$$

The auxiliary function in (52) has been written as $\left( \frac{x-x_0}{|\bar{r}-\bar{r}_0|^4} \right) e^{-(w-w_0)\bar{\nabla}}$ to emphasise its status as right-regular, but it is identical to $e^{-(w-w_0)\bar{\nabla}} \left( \frac{x-x_0}{|\bar{r}-\bar{r}_0|^4} \right)$ due to the generating function being scalar and hence the explicit transcendental expression (47) applies.

## Extension to Biquaternions

Moving now to biquaternions, a spacetime location is specified by an Hermitian biquaternion,

$$q = t + i(xI + yJ + zK) \tag{54}$$

where we use the symbol $t$ for the scalar part as we are now dealing with a natural norm which is Minkowskian, and hence $t$ can be interpreted as physical time.

## Biregular Functions

For biquaternionic functions we introduce the property of biregularity. For this purpose we define the Hermitian operator $D_H = \partial_t + i\bar{\nabla}$, where the subscript $H$ denotes "Hermitian", because this operator is such that $D_H^{\#} = D_H^* = \partial_t - i\bar{\nabla}$. We define a biquaternion-valued function, $f(t, x, y, z)$, to be biregular if,

$$D_H f(t, x, y, z) = 0 \tag{55}$$

Hence, whereas regular functions obey Laplace's equation in $\mathbb{R}^4$, $\nabla_4^2 f = 0$, biregular functions obey the wave equation,

$$D_H^{\#} D_H f = (\partial_t^2 - \nabla_3^2) f = 0 \tag{56}$$

Similarly we can define conjugate biregular functions by $D_H^{\#} f = 0$, right-biregular functions by $f \overleftarrow{D}_H = 0$ and conjugate right-biregular functions by $f \overleftarrow{D}_H^{\#} = 0$, all of which also obey the wave equation, (56).

Analogous to (42), biregular functions can be written in terms of a generating function as,

$$f(t, x, y, z) = e^{-it\bar{\nabla}} f(0, x, y, z) \tag{57}$$

Similarly conjugate biregular functions can be written $f(t, x, y, z) = e^{it\bar{\nabla}} f(0, x, y, z)$, whilst right-biregular functions can be written $f(t, x, y, z) = f(0, x, y, z) e^{-it\overleftarrow{\bar{\nabla}}}$ and conjugate right-biregular functions as $f(t, x, y, z) = f(0, x, y, z) e^{it\overleftarrow{\bar{\nabla}}}$.

Just as for regular functions, the theorem follows that if $g$ and $f$ are both biregular and $g$ has a scalar generating function, then $gf$ is also biregular, the proof being as (43,44).

## Biquaternionic Forms

The natural 3-surface element, or 3-form, is now the Hermitian biquaternion analogue of (22),

$$S_H = dx \wedge dy \wedge dz - i(I dw \wedge dy \wedge dz + J dw \wedge dz \wedge dx + K dw \wedge dx \wedge dy) \tag{58}$$

whilst the 4-form, or volume element in spacetime, is $V_H = dt \wedge dx \wedge dy \wedge dz$, which is formally the same as $V_q$ as the latter is a real scalar and hence Hermitian. The key result that permits Cauchy-like surface-independent integrals over regular functions to be demonstrated from the Stokes-Cartan Theorem, (19), is (30). The analogue holds for biregular functions, noting again that $dS_H = 0$ so we get,

$$d(S_H f) = -V_H D_H f \tag{59}$$

which is therefore zero for biregular functions. Similarly, for conjugate-regular functions, $d(S_H^\# f) = -V_H D_H^\# f = 0$; for left-regular functions $d(f S_H) = (f \overleftarrow{D}_H) V_H = 0$, and for conjugate left-regular functions $d(f S_H^\#) = (f \overleftarrow{D}_H^\#) V_H = 0$.

Analogous to the case of regular quaternion functions, the condition (55) that a biquaternion function is biregular is equivalent to the existence of a "derivative", $g$, such that,

$$\frac{1}{2} d(dq \wedge dq f) = S_H g \tag{60}$$

$\frac{1}{2} dq \wedge dq$ is again given by (23) except for a change of sign. The demonstration that (60) is equivalent to (55) follows from $\frac{1}{2} d(dq \wedge dq f)$ evaluating to $-i S_H \partial_t f$ if (55) holds, analogous to (24,25).

## Cauchy-like Integral Theorems for Biregular Biquaternionic Functions

It follows that,

$$\oiiint_{\partial V} S_H f(t, x, y, z) = 0 \tag{61}$$

for any function which is biregular within and on the closed boundary $\partial V$. To emphasise: this follows because the Stokes-Cartan Theorem, (19), shows that this integral is equal to the integral over $V$ of $d(S_H f) = -V_H D_H f = 0$.

We now seek an extension of (48) into the biquaternion domain. Since we require a biregular auxiliary function, in view of (57) it is natural to explore,

$$\tilde{g} = e^{-it\bar{\nabla}} \left( \frac{x}{r^4} \right) \tag{62}$$

as the analogue of (46). As would be expected from the replacement $w \to it$, the equivalent of (47) is,

$$\tilde{g}(t, x, y, z) = \frac{x}{(t^2 - r^2)^2} - iI \frac{t}{2r^2} \left\{ \frac{1}{r^2 - t^2} + \frac{1}{rt} \tanh^{-1}\left(\frac{t}{r}\right) \right\} + 4ixt \frac{\bar{r}}{r^4} \left\{ \frac{5r^2 - 3t^2}{8(t^2 - r^2)^2} + \frac{3}{8rt} \tanh^{-1}\left(\frac{t}{r}\right) \right\} \tag{63}$$

It follows that, as the generating function of $\tilde{g}$ is scalar, that $\tilde{g} f$ is biregular when $f$ is biregular, except at the singularity of $\tilde{g}$. Hence the integral,

$$\oiiint_{\partial V} S_H \left[ e^{-i(t-t_0)\bar{\nabla}} \left( \frac{x - x_0}{|\bar{r} - \bar{r}_0|^4} \right) \right] f(t, x, y, z) = \xi f(t_0, x_0, y_0, z_0) \tag{64}$$

over the boundary of a simply connected region, $V$, is independent of $V$ provided that the singular point, $q_0$, is fully contained within it, and $f$ is biregular everywhere within $V$ and on its boundary.

However, we have yet to evaluate the constant $\xi$ or even to demonstrate that it is finite and non-zero. There is a problem here that did not arise in the quaternionic case, namely that $\tilde{g}$ is not singular only at a point but everywhere on the light cone, i.e., $t^2 - r^2 = 0$. Because of this difficulty we spell out in detail how the integral,

$$\xi = \oiiint_{\partial V} S_H e^{-i(t-t_0)\bar{v}} \left(\frac{x-x_0}{|\bar{r}-\bar{r}_0|^4}\right) \tag{65}$$

is evaluated. As it stands, (65) is ambiguous. It is made unique, and finite, only by acknowledging that it must be interpreted as a contour integral in the Argand plane and the path of the contour around the singularities defined.

## Evaluation of (65) on a Long Prismatic Surface

The origin can be shifted to $q_0$ without changing the integral. As (65) is independent of the integration surface we opt to evaluate it on a prismatic surface defined as follows.

- Consider a closed 2-surface $\partial V_3$ enclosing a 3D region $V_3$ of the spatial part of biquaternion space, i.e., the 3-vector part. This region $V_3$ is extruded parallel to the time axis from $-t_1$ to $+t_1$, with the intention of considering $t_1 \to \infty$.

- We can take $\partial V_3$ to be a spherical 2-surface of finite radius $\rho$, centred on $q_0 = 0$. The 3-surface element $S_H$ on the curved cylindrical surface is then the product of this $\partial V_3$ with $dt$, which is therefore oriented parallel to the spatial radial direction.

- The biquaternion integration measure $S_H$ on the cylindrical surface can thus be represented in spherical polar coordinates as,

$$S_H \to i(I\sin\vartheta\cos\varphi + J\sin\theta\sin\varphi + K\cos\theta)dt r^2 \sin\theta d\theta d\varphi$$

being careful to note that we will need to integrate along the time coordinate as well as over the spherical 2-surface.

- The cylindrical surface is closed by adding 'end caps' at times $\pm t_1$. These 'end caps' are simply the spatial 3-volumes $V_3$ within radius $\rho$ at times $\pm t_1$. At $+t_1$ the outward normal is positive, whereas at $-t_1$ the outward normal is negative. Hence the integration measure, $S_H$, becomes simply $dxdydz$ at $+t_1$ but $-dxdydz$ at $-t_1$.

We now evaluate the integral (65) where the integrand is given by (63), taking the two parts in turn.

### End Caps

The end caps lie at $t = \pm t_1$ and we are interested in the limit $t_1 \to \infty$. The first term in (63), the scalar component, is thus zero on the end caps. The same is true for the first terms within each of the two $\{...\}$. In the remaining terms the $\tanh^{-1}\left(\frac{t}{r}\right)$ becomes $\tanh^{-1}(\pm\infty)$. As the integrand is equal and opposite on the two end caps, the two ends do *not* cancel but add.

Other than that factor of $\tanh^{-1}(\pm\infty)$ that leaves us with the requirement to integrate over the volume of a sphere of radius $\rho$ the function,

$$-\frac{iI}{2r^3} + \frac{3ix\bar{r}}{2r^5}$$

The first term integrates to $-2i\pi I \int_0^\rho \frac{dr}{r}$ which is divergent. The second term integrates as follows,

$$\frac{3i}{2} \int \frac{1}{r^5} r\sin\theta\cos\varphi . r(I\sin\vartheta\cos\varphi + J\sin\theta\sin\varphi + K\cos\theta) r^2 dr . \sin\theta d\theta d\varphi$$

The $\varphi$ integral kills the $J$ and $K$ terms and we are left with,

$$\frac{3i}{2} I \int_0^\rho \frac{dr}{r} \int_0^\pi \sin^3\theta d\theta \int_0^{2\pi} \cos^2\varphi d\varphi = \frac{3i}{2} I \cdot \frac{4}{3} \cdot \pi \int_0^\rho \frac{dr}{r} = 2i\pi I \int_0^\rho \frac{dr}{r}$$

which cancels with the first term, above. So the end caps integrate to zero.

Cylindrical Surface

Considering firstly the scalar term in the integrand, (63), we need to evaluate,

$$\int i(I\sin\vartheta\cos\varphi + J\sin\theta\sin\varphi + K\cos\theta) \frac{r\sin\theta\cos\varphi}{(t^2-r^2)^2} dt\, r^2\sin\theta\, d\theta d\varphi$$

The $\varphi$ integral kills the $J$ and $K$ terms and we are left with,

$$iI \int_0^\pi \sin^3\theta d\theta \int_0^{2\pi} \cos^2\varphi d\varphi \int_{-\infty}^{+\infty} \frac{r^3}{(t^2-r^2)^2} dt = \frac{4\pi}{3} iI \int_{-\infty}^{+\infty} \frac{r^3}{(t^2-r^2)^2} dt \qquad (66)$$

We will come back to the evaluation of the integral $\int_{-\infty}^{+\infty} \frac{r^3}{(t^2-r^2)^2} dt$ later as this involves careful handling to avoid being divergent due to the singularities that now occur at $t = \pm r$. We did not have this difficulty when working in Euclidean quaternion space.

We now show that the two vector parts of (63) both integrate to zero.

For the first of the vector terms this is immediately clear because it contains no angular dependence, and so inevitably integrates to zero over the vectorial integration measure oriented parallel to $\hat{r}$.

For the second vector term, its product with the integration measure involves $\hat{r}^2 = -1$, which eliminates the angular dependence of those factors. The remaining angular-dependent terms are $x.\sin\theta d\theta d\varphi = r\sin^2\theta\cos\varphi d\theta d\varphi$ which is killed by the $\varphi$ integral.

Hence, the required integral reduces to just (66). We now need to face the issue of the singularity in that integral over $t$.

**Evaluation of (66)**

The indefinite integral is,

$$\int \frac{r^3}{(t^2-r^2)^2} dt = \int \frac{d\tau}{(\tau^2-1)^2} = \frac{\tau}{2(1-\tau^2)} + \frac{1}{4} \log\left(\frac{1+\tau}{1-\tau}\right) \qquad \tau = \frac{t}{r} \qquad (67)$$

However, the integrand is singular at $\tau = \pm 1$. If we adopt the principal value,

$$\mathbb{P} \int_{-\infty}^\infty \frac{d\tau}{(\tau^2-1)^2} = 2\,\mathbb{P} \int_0^\infty \frac{d\tau}{(\tau^2-1)^2} = 2 Lim_{\varepsilon \to 0} \left( \int_0^{1-\varepsilon} \frac{d\tau}{(\tau^2-1)^2} + \int_{1+\varepsilon}^\infty \frac{d\tau}{(\tau^2-1)^2} \right) \qquad (68)$$

then we find this is divergent as $1/\varepsilon$. However, since this is part of a boundary which is required to be closed, the only option is to skirt around the singularities by interpreting them as poles in the Argand plane of a complex function.

Hence we reinterpret the integral to be,

$$\int_{-\infty}^{\infty} \frac{d\tau}{(\tau^2-1)^2} \to \oint_C \frac{dz}{(z^2-1)^2} \tag{69}$$

where it is understood that the contour is closed via a semicircle pushed off to infinity in the upper half plane. However, there is still an ambiguity in how the integration contour detours around the two poles, at $z = \pm 1$. We choose to define the contour as including the pole at $z = +1$ but excluding the pole at $z = -1$, as shown in Figure 1.

**Figure 1: Integration contour in the Argand plane for (69)**

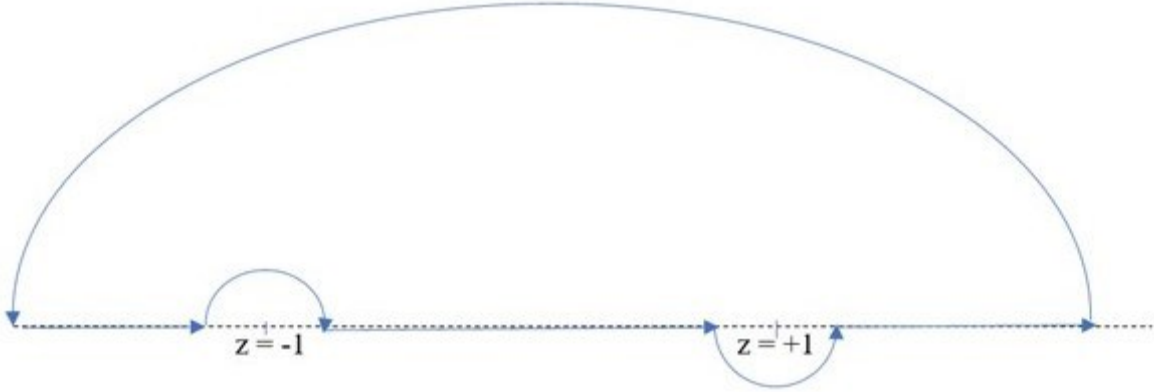

where it is understood that the upper half-circle is pushed to infinity, and hence contributes nothing to the contour integral. The contour integral in this limit does, therefore, become the desired real integral – but subject to our choice of which pole to include within the contour and which to exclude. If either both or neither of the poles are included, the integral would be zero as the two poles have equal and opposite residues. If only one pole is included, then which is chosen determines the sign of the result but not its magnitude.

The poles are of second order. Hence the residue at $z = 1$ is given by,

$$\text{Residue: } \frac{d}{dz}\left[(z-1)^2 \frac{1}{[(z-1)(z+1)]^2}\right]\bigg|_{z=1} = -2(z+1)^{-3}\big|_{z=1} = -\frac{1}{4} \tag{70}$$

The contour integral is thus $2\pi i$ times this residue, hence $-i\pi/2$. (66) thus gives the value of the integral (65) to be $\xi = \frac{2}{3}\pi^2 I$, exactly the same as the quaternion case, (48). Hence, the integral theorem for biregular functions is,

$$\oiiint_{\partial V} S_H\left[e^{-i(t-t_0)\bar{\nabla}}\left(\frac{x-x_0}{|\bar{r}-\bar{r}_0|^4}\right)\right] f(t,x,y,z) = \frac{2}{3}\pi^2 I f(t_0, x_0, y_0, z_0) \tag{71}$$

Similar equations apply when $x$ and $I$ are replaced by $y, z$ and $J, K$, of course. In the same manner we find the corresponding integral theorems for conjugate biregular, right-biregular and conjugate right-biregular functions analogous to (51-53) by replacing $(w - w_0)$ in the exponents with $i(t - t_0)$. As one example, for conjugate right-biregular functions we have,

$$\oint_{\partial V} f(q)\left[\left(\frac{x-x_0}{|\bar{r}-\bar{r}_0|^4}\right) e^{i(t-t_0)\bar{\nabla}}\right] S_q^\# = \oint_{\partial V} f(q)\left[e^{i(t-t_0)\bar{\nabla}}\left(\frac{x-x_0}{|\bar{r}-\bar{r}_0|^4}\right)\right] S_q^\# = \frac{2\pi^2}{3} f(q_0) I \tag{72}$$

Fueter's Theorem also extends to biregular functions. The auxiliary function is now,

$$-e^{-it\bar{\nabla}}\left(\frac{\bar{r}}{r^4}\right) = \frac{iq}{|q|^4} = \frac{i(t+i\bar{r})}{(t^2-r^2)^2} \tag{73}$$

So Fueter's Theorem for a biregular function is,

$$\oiiint_{\partial V} \frac{i(q-q_0)}{|(q-q_0)|^4} S_H f(t,x,y,z) = 2\pi^2 I f(t_0, x_0, y_0, z_0) \tag{74}$$

This works because the analogue of (34) continues to apply, i.e.,

$$d(g S_H f) = (\widetilde{D}_H g) V_H f + g V_H D_H f \tag{75}$$

So that, if $g$ is a right-biregular function and $f$ is a left-biregular function, then the integral $\oiint_{\partial V} g(q) S_H f(q)$ is surface-independent and zero if singularities are excluded. As before, the auxiliary function, (73), is both left and right biregular, so the theorem (74) applies subject only to confirming the value of the constant factor. The required integral is,

$$\oiiint_{\partial V} \frac{i(t+i\bar{r})}{(t^2-r^2)^2} i\hat{r} r^2 d\Omega dt = \oiiint_{\partial V} \frac{(-t\hat{r}+ir)}{(t^2-r^2)^2} r^2 d\Omega dt \tag{76}$$

Opting to evaluate on an infinitely long prismatic surface formed by extruding a spatial 2-sphere along the time axis, as before, the integrand is zero on the end caps whilst the $\hat{r}$ term integrates to zero over the 2-sphere. This leaves $4\pi i \int_{-\infty}^{\infty} \frac{d\tau}{(\tau^2-1)^2}$. This reduces to the same integral as above, and adopting the same interpretation, $\int_{-\infty}^{\infty} \frac{d\tau}{(\tau^2-1)^2} \to \oint_C \frac{dz}{(z^2-1)^2} = -i\frac{\pi}{2}$ gives the final result for (76) to be $2\pi^2$, confirming (74).

## Evaluation of (65) on a Wide Prismatic Surface

To confirm the surface independence of (65) by direct calculation on an alternative surface we again choose a prismatic surface with a 2-sphere spatial boundary, but now of finite extent in time (i.e., finite $\pm t_1$) but of infinite spatial extent, i.e., $\rho \to \infty$.

On the cylinder boundary, as $\rho \to \infty$, the $tanh^{-1}$ terms in the integrand of (65), i.e., (63), become zero and the remaining terms are all of order $1/\rho^4$, and hence tend to zero when multiplied by the integration measure $\rho^2 dt d\Omega$.

The scalar term in (63) is even in $t$ and hence cancels between the two end caps due to their opposite orientation. The remaining terms are odd in $t$ and hence contribute equally on the two ends caps.

Carrying out the angular integrations shows that the two $tanh^{-1}$ terms cancel, as follows. The factor multiplying the first of these terms after integrating over $d\Omega$ is, for the end cap at positive $t$,

$$-iI\frac{4\pi}{2r^3} \tag{77}$$

For the second of the two $tanh^{-1}$ terms, its factor after angular integration on the same end cap is,

$$\int_{\varphi=0}^{2\pi} \int_{\theta=0}^{\pi} \frac{3i}{2r^3} sin\theta cos\varphi (I sin\theta cos\varphi + J sin\theta sin\varphi + K cos\theta) sin\theta d\theta d\varphi \tag{78}$$

The $\varphi$-integral kills the $J$ and $K$ terms. The integral over $sin^3\theta$ produces the factor $4/3$ whilst the $cos^2\varphi$-integral gives a factor of $\pi$. Hence (78) cancels with (77) and we are spared having to integrate over the awkward $tanh^{-1}$ terms.

What remains to be evaluated, after carrying out the angular integrations, is,

$$2\int_0^\infty r^2 dr \left(-iI\frac{2\pi t}{r^2(r^2-t^2)} + iI\frac{2\pi(5r^2-3t^2)t}{3r^2(r^2-t^2)^2}\right) = \int_0^\infty dz \left(-iI\frac{4\pi}{(z^2-1)} + iI\frac{4\pi(5z^2-3)}{3(z^2-1)^2}\right)$$

where the factor of 2 accounts for both end caps and $z = r/t$. As before we avoid the singularities in the integrand by suitable choice of contour in the Argand plane, so the desired integral becomes (noting that the integrand is even in $z$),

$$\frac{1}{2}\oint_{C'} \left(-iI\frac{4\pi}{(z^2-1)} + iI\frac{4\pi(5z^2-3)}{3(z^2-1)^2}\right) dz$$

Here the whole of the real line is included in contour $C'$, except for skirting around the poles at $z = \pm 1$, and the contour is closed in the upper half plane noting that the infinite semicircle contributes nothing. However, to be consistent we must now include the poles at $z = -1$ and exclude those at $z = +1$. This is because the previous integration was over integration variable $z = t/r$ whereas now the integration variable is $z = r/t$, and a small positive imaginary contribution to the former becomes a small negative imaginary part in the latter, and *vice-versa*. This reverses which pole is included in the contour.

The residue of $\frac{1}{(z^2-1)}$ at $z = -1$ is $-\frac{1}{2}$.

The residue of $\frac{(5z^2-3)}{(z^2-1)^2}$ at $z = -1$ is $\frac{d}{dz}\left((z+1)^2 \frac{(5z^2-3)}{(z^2-1)^2}\right)\bigg|_{z=-1} = -2$

Hence, we get that (65) evaluates to,

$$\xi = \frac{1}{2} \times 2\pi i \left[-4\pi iI \times -\frac{1}{2} + \frac{4\pi}{3} iI \times -2\right] = \frac{2}{3}\pi^2 I$$

Thus confirming the previous result and hence the integral theorem (71) for biregular functions.

## Conclusion

Novel integral theorems of Cauchy type have been presented for quaternion functions which are regular (or right-regular or their conjugates) in (48 – 53). These theorems appear to be distinct from Fueter's Theorem, (32), in that they rely on different features in order to be integration surface independent and are not obviously inter-convertible.

Corresponding theorems have been shown to apply for biquaternion functions which are biregular (or right biregular or their conjugates) in (71, 72), but these integrals require careful specification of the integration contours in the Argand plane to include just one of the two poles.

The equivalent of Fueter's Theorem also applies for biregular biquaternion functions, (74), providing the same caution is exercised regarding the integration contour deployed around the poles of the integrand.